\documentclass[12pt]{amsart}

\usepackage{amssymb}

\newtheorem{theorem}{Theorem}[section]
\newtheorem{proposition}[theorem]{Proposition}
\newtheorem{lemma}[theorem]{Lemma}
\newtheorem{remark}[theorem]{Remark}

\newtheorem{definition}[theorem]{Definition}

\numberwithin{equation}{section}

\begin{document}
\baselineskip=15.5pt

\title[Moduli of Vortices and Grassmann Manifolds]{Moduli of
Vortices and Grassmann Manifolds}

\author[I. Biswas]{Indranil Biswas}

\address{School of Mathematics, Tata Institute of Fundamental
Research, Homi Bhabha Road, Bombay 400005, India}

\email{indranil@math.tifr.res.in}

\author[N. M. Rom\~ao]{Nuno M. Rom\~ao}

\address{Centre for Quantum Geometry of Moduli Spaces,
Institut for Matematiske Fag,
Aarhus Universitet, Ny Munkegade bygn.\ 1530,
DK-8000 \AA rhus C, Denmark}

\email{nromao@imf.au.dk}



\date{}

\begin{abstract}

We use the framework of Quot schemes to
give a novel description of the moduli spaces of stable $n$-pairs, also interpreted as 
gauged vortices on a
closed Riemann surface $\Sigma$ with target 
${\rm Mat}_{r\times n}(\mathbb{C})$, where $n\ge r$. We then show that
these moduli spaces embed canonically into certain Grassmann manifolds, and
thus obtain natural K\"ahler metrics of Fubini--Study type. These spaces are smooth
at least in the local case $r=n$. For abelian local vortices we prove that, if a 
certain ``quantization'' condition
is satisfied,  the embedding can be chosen in such a way that the induced
Fubini--Study structure realizes the K\"ahler class
of the usual $L^2$ metric of gauged vortices.

\end{abstract}

\maketitle

\section{Introduction}\label{intro}

Gauged vortices are configurations of static, stable fields arising in various
classical field theories on a Riemann surface $\Sigma$. These objects were first
studied as topological solitons of the abelian Higgs model, for which
vortex solutions have a distinctive particle-like behavior --- they are
labelled by divisors on $\Sigma$, which specify the precise
locations of the cores of individual objects superposing nonlinearly to yield each vortex 
configuration~\cite{JafTau}.
In this setting, there is typically a moduli space of all vortices with a given
topology, modelled on the space of effective divisors with a fixed degree. This is a
smooth manifold endowed with a
complex structure induced from the one specified on $\Sigma$.
More recently, models for vortices with nontrivial internal structure have been
considered, but in the various generalizations it has remained a challenge to
understand the corresponding moduli spaces in a satisfactory way.

We shall focus on vortices on a closed Riemann surface $\Sigma$ with
target (or internal) space consisting of the vector space 
${\rm Mat}_{r \times n}(\mathbb{C})$ of complex
$r \times n$ matrices, where $n \ge r$. These have been called nonabelian
vortices in the literature, even though the special situation $r=1$ corresponds
to an abelian gauge theory. If $n>r$, one sometimes speaks of semilocal
vortices, whereas $n=r$ is known as the local case.
The geometric framework is as follows. Let $e^2$ be any positive real number. Assume that
we fix a K\"ahler form $\omega_\Sigma$ on $\Sigma$, as well as a Hermitian
metric on a complex vector bundle $E \longrightarrow \Sigma$ of rank 
$r$. A vortex
is a pair $(A,\phi)$ consisting of a unitary connection $A$ on the bundle, together
with a section $\phi$ of the direct sum $E^{\oplus n} \longrightarrow 
\Sigma$, satisfying the vortex equations
\begin{eqnarray}
&& \bar\partial_A \phi =0, \label{vor1} \\
&& * F_{A} + e^2 \mu \circ \phi = 0. \label{vor2}
\end{eqnarray}
Here, $\bar \partial_A$ denotes the holomorphic structure~\cite{DonKro} on $E^{\oplus n}$
defined by the
connection $A$ and the complex structure on $\Sigma$, $*$ is the Hodge star
of the K\"ahler metric associated to $\omega_\Sigma$, $F_A:={\rm d}A+\frac{1}{2}[A,A]$ is the curvature
of $A$ and $\mu$ denotes a moment map
\[
\mu: {\rm Mat}_{r \times n}(\mathbb{C}) \longrightarrow \mathfrak{u}(r)^{*} \cong \mathfrak{u}(r)
\]
of the Hamiltonian action of the reduced structure group ${\rm U}(r)$ on the fibers
of $E^{\oplus n} \longrightarrow \Sigma$ by multiplication on the left. 
We use the
Killing form on $\mathfrak{u}(r)$ to identify the Lie algebra with its dual.
Notice that $\mu$ is specified only up to addition of scalar matrices, and
following standard conventions we shall write
\[
\mu(w) = -\frac{\sqrt{-1}}{2}(w w^{\dag} - \tau I_r),
\]
where $\tau$ is a fixed real number and $I_r$ is the $r \times r$ identity matrix. 

The vortex equations \eqref{vor1}--\eqref{vor2} first appeared in the 
work~\cite{BDW} of Bertram, Daskalo\-poulos
and Wentworth computing the Gromov--Witten
invariants of Grassmannians: the moduli space of holomorphic maps from 
a compact Riemann surface to a Grassmannian embeds into the moduli space of stable holomorphic $n$-pairs. The latter can be identified with the space of gauge-equivalence classes of solutions to the vortex equations
above, under suitable stability criteria depending on the parameter $\tau$ and
the topology. This is an example of what is generically known as the Hitchin--Kobayashi correspondence,
which goes back to~\cite{UhlYau}. 
Among other things, the authors of~\cite{BDW} described how the moduli 
space of vortices changes birationally when the parameter $\tau$ crosses 
certain critical values, a phenomenon familiar from earlier work of Thaddeus on moduli of stable pairs~\cite{Tha}.
There is also a useful description of the moduli spaces via infinite-dimensional
symplectic reduction (in the spirit of~\cite{AtiBot}), which naturally produces a
K\"ahler structure from the $L^2$ inner product on the space of fields; 
for abelian vortices, this was described by Garc\'ia-Prada 
in~\cite{GaDEP}.
By now, a whole body of rather well-established technology that
reproduces results of this type has been developed for objects that
are analogous to vortices on the gauge-theory side of the Hitchin--Kobayashi correspondence, and $n$-pairs on the other side. The objects on the algebraic-geometric side are often referred to by the name of augmented bundles, of which Higgs bundles and coherent systems are other important examples; we refer the reader to~\cite{BDGW} for a clear overview.

Physicists have also been interested in the generalized vortex equations
(\ref{vor1})--(\ref{vor2}). Their solutions realize certain configurations of branes in string theory on the one hand, and also feature in models for confinement in
QCD~\cite{EINOS,TonQVS}. Here, one focus of interest is to obtain descriptions of the
moduli spaces as explicit as possible, including concrete parametrizations, as well
as to understand natural Hamiltonian systems on the moduli spaces or their cotangent bundles. Much of the work done 
assumes $\Sigma=\mathbb{C}$, for which there is nothing like a Hitchin--Kobayashi
correspondence, but alternative constructions have been proposed which
rely on certain mathematical conjectures. More recently, Baptista presented a
rigorous description of the moduli space of local vortices when $\Sigma$ is
compact, describing a stratification of the moduli spaces in terms of spaces of internal 
structures~\cite{BatNAV}. From Baptista's description, holomorphic 
matrices representing vortex solutions up to unitary gauge 
transformations can be readily constructed. From our perspective, his 
work has the slight disadvantage of depending on auxiliary
structure, namely the choice of an inner product on $\mathbb{C}^n$, and it is
also difficult to see how the different strata are patched together.

In this paper, we make use of the Hitchin--Kobayashi correspondence 
of~\cite{BDW} to
describe moduli spaces of solutions to the vortex equations \eqref{vor1}--\eqref{vor2} modulo gauge equivalence,
\begin{equation}\label{moduli}
{\mathcal M}_\Sigma \,= \,{\mathcal M}_\Sigma (n,r,d),
\end{equation}
where $d={\rm deg} (E)$ is the degree of $E \longrightarrow \Sigma$, in 
terms of certain Quot schemes parametrizing
holomorphic $n$-pairs.
The idea of Quot (or quotient) schemes goes back to Grothendieck~\cite{Gr}
and has had numerous applications to moduli problems.
Given a coherent sheaf and a polynomial, the Quot scheme is a projective scheme
of finite type that parametrizes all quotients of the given sheaf for which the
Hilbert polynomial~\cite{EisHar} is the given polynomial.

Starting with an ample line bundle ${\mathcal L} \longrightarrow 
\Sigma$, we shall 
produce a holomorphic embedding of the moduli space ${\mathcal M}_\Sigma$ into a Grassmann manifold;
it follows that ${\mathcal M}_\Sigma$ is projective.
A Hermitian structure on ${\mathcal L}$ then induces a K\"ahler metric of Fubini--Study type on the
moduli space.
The perspective of Quot schemes has the advantage of being global in nature, and also
well-suited to address general questions such as smoothness. 
We shall also see how it allows a straightforward calculation of the dimension.
These properties can also be recovered from more general results 
scattered in the literature on Gromov--Witten invariants~\cite{OkoTel, 
BDW}.

The simplest example of our class of embeddings into Grassmann manifolds occurs when
we set $n=r=1$; more background on the geometry of the moduli space of vortices in this well-studied case shall be given in Section~\ref{bgabellocal} 
below. Then we have
${\mathcal M}_{\Sigma} \cong {\rm Sym}^d(\Sigma)$, where
$d\,=\,{\rm deg}(E)$ is the vortex number~\cite{Br}.
In this setting, one might hope that a suitable
choice of hermitian metric on ${\mathcal L}$ will induce a Fubini--Study metric
which is related to the usual $L^2$ metric on the moduli space of vortices.
We shall show that, if a certain quantization condition holds, then it turns out that
the two corresponding K\"ahler classes are cohomologous; this is the content of
our Theorem~\ref{cohomL2} below. The K\"ahler class of the
moduli space of local abelian vortices was calculated in 
\cite{ManNasVVMS}.

\section{Stability and the Hitchin--Kobayashi correspondence}\label{stabHitKob}

Let $\Sigma$ be a compact connected Riemann surface of genus $g$. Fix
a K\"ahler form $\omega_\Sigma$ on $\Sigma$, so $\omega_\Sigma$ is
a positive $(1,1)$-form; it is automatically closed. We will denote by
\begin{equation}\label{vo.}
{\rm Vol} (\Sigma)\,:=\, \int_{\Sigma}\omega_\Sigma
\end{equation}
the total area of the surface determined by $\omega_\Sigma$.

We briefly sketch the results in \cite{BDW} establishing the
Hitchin--Kobayashi correspondence between solutions of 
\eqref{vor1}--\eqref{vor2} up to gauge transformations, and stable 
$n$-pairs $(E,s)$ up to isomorphism. We begin by recalling the following

\begin{definition}
An {\em $n$-pair} on the Riemann surface $\Sigma$ is a pair of the form $(E,s)$, where
$E \longrightarrow \Sigma$ is a holomorphic vector bundle  and 
$s\,\in\, H^0(\Sigma,\,E^{\oplus n})$. Two $n$-pairs $(E,s)$ and $(E',s')$ are said to be 
{\em isomorphic} if
there is an isomorphism of holomorphic vector bundles $\psi: E \longrightarrow 
E'$ over $\Sigma$ such that $\psi^*s'=s$.
\end{definition}
In this paper, we will denote by $r\,=\,{\rm rk}(E)$ the rank of a fixed class of 
vector bundles $E$ over
$\Sigma$, when no confusion will arise.

The basic mechanism of the correspondence is modeled on Donaldson's famous
proof of the Narasimhan--Seshadri theorem~\cite{DonNS}. Suppose that we are given
an $n$--pair $(E,s)$. A holomorphic vector bundle $E \longrightarrow 
\Sigma$ with a Hermitian structure has a unique connection $A$ 
preserving the Hermitian structure whose $(0,1)$-part coincides with the Dolbeault
operator defining the holomorphic structure; this connection $A$ is 
known as the \textit{Chern connection} \cite[pp.~191--192, Proposition 
5]{At}, \cite[p.\ 273]{CCL}. For a $C^\infty$ section $\phi$ of
$E^{\oplus n} \longrightarrow \Sigma$, the pair $(A,\phi)$ is a solution of (\ref{vor1})
if and only if $\phi$ is holomorphic. So we start by taking
$\phi= s\,\in\, H^0(\Sigma,\,E^{\oplus n})$.
Complex gauge transformations preserve equation (\ref{vor1}), and one can ask whether the  complex gauge orbit through
this initial pair $(A,\phi=s)$ contains a solution of equation (\ref{vor2}), which
itself is invariant only under unitary gauge transformations. The answer is that
this occurs if and only if the $n$--pair $(E,s)$ is $\tau$--stable in a sense that we
will explain shortly, for the value of $\tau$ appearing in equation (\ref{vor2}). This solution is unique up to unitary gauge transformations, and therefore we obtain an injective map from the
moduli space of $\tau$-stable $n$--pairs to the moduli space of vortices. Conversely,
a vortex $(A,\phi)$ in this geometric setting determines an $n$--pair: $E$ is the bundle where each
component of $\phi$ takes values, with holomorphic structure 
on $E$ determined by the connection
$A$ and the complex structure on $\Sigma$. Clearly, one obtains isomorphic
$n$-pairs $(E,\phi)$ when the original vortex $(A,\phi)$ undergoes unitary
gauge transformations, and one can check that they are still $\tau$--stable.

The stability condition that is appropriate to relate $n$-pairs and vortices
was spelled out in~\cite{BDW,BDGW}, using the analysis for stable pairs in~\cite{Br2}. Fixing $\tau$, one says that an $n$-pair $(E,\phi)$ is
$\tau$-stable if the following two conditions hold:
\begin{itemize}
\item[(i)]
$4 \pi\,{\deg}(E')/{\rm rk}(E')\,<\, {\tau e^2 \,{\rm Vol}(\Sigma)} \;$  
for all holomorphic subbundles $E' \subseteq E$, and
\item[(ii)]
$4 \pi\, {\deg}(E/E_s)/ {\rm rk}(E/E_s)\,
>\, {\tau e^2 \,{\rm Vol}(\Sigma)} \;$ for all holomorphic subbundles 
$E_s \subsetneq E$ containing
all the component sections of ${s}$.
\end{itemize}
(${\rm Vol} (\Sigma)$ is defined in \eqref{vo.}; unlike~\cite{BDW,BDGW}, 
we do not require
this area to be normalized.) Notice that, when $E'=E$, 
condition (i) is necessary for vortex solutions to exist for a given $\tau$: this
follows from integrating equation (\ref{vor2}) over $\Sigma$.

Now suppose that $n \ge r = {\rm rk}(E)$, and that $\phi \in H^0(\Sigma,\,E^{\oplus n})$ has maximal rank generically on
$\Sigma$. Then there is no proper subbundle of $E$ containing all the components
of $\phi$, and the second condition above is empty. Going through the
argument in the proof of Proposition~3.14 in~\cite{BDW}, one can show 
that, under the same assumptions, the inequality
\begin{equation}\label{stable}
\tau e^2\, {\rm Vol}(\Sigma)  > 4 \pi \, {{\deg} (E)}
\end{equation} 
is equivalent to the first condition for $\tau$-stability. 
Throughout this paper, when the
topology of $E \longrightarrow \Sigma$ has been fixed, as well as a 
K\"ahler
structure on $\Sigma$, we shall only deal with the vortex equations
\eqref{vor1}--\eqref{vor2} with values of $\tau$ satisfying (\ref{stable}).
Then we can focus purely on $n$-pairs and their algebraic geometry to
describe the moduli spaces in (\ref{moduli}).

\section{Holomorphic sections of a direct sum}\label{sec1}

We shall from now on take the algebraic-geometric point of view on the moduli space
of vortices provided by the Hitchin-Kobayashi correspondence explained in Section~\ref{stabHitKob}.

As before, let $E\, \longrightarrow\, \Sigma$ be a holomorphic vector bundle
of rank $r$. Choose an integer $n\, \geq\, r$. Let
\begin{equation}\label{e1}
s\, \in\, H^0(\Sigma,\, E^{\oplus n})\,\cong \, H^0(\Sigma,\, E)^{\oplus n}
\end{equation}
be a holomorphic section. Let $s_i\,\in\,
H^0(\Sigma,\, E)$ be the image of $s$ for the projection
$E^{\oplus n}\, \longrightarrow\,E$ to the $i$--th factor.
So, $s\, =\, (s_1\, ,\ldots\, , s_n)$. Let
\begin{equation}\label{e2}
f_s\, :\, {\mathcal O}^{\oplus n}_\Sigma\,\longrightarrow\, E
\end{equation}
be the homomorphism defined by $(x\, ; c_1\, ,\cdots\, ,c_n)
\, \longmapsto\, \sum_{i=1}^n c_i\cdot s_i(x)$, where $x\, \in\,
\Sigma$ and $c_i\, \in\, \mathbb C$. The image ${\rm im}(f_s)$
is locally free; however, $\text{im}(f_s)$ need not be a subbundle of $E$.

\begin{definition}\label{def1}
{\rm Let}
$$
H^0(\Sigma,\, E^{\oplus n})_0\,\subset\, H^0(\Sigma,\, E^{\oplus n})
$$
{\rm be the subset consisting
of sections $s$ as in \eqref{e1} such that
the rank of the vector bundle ${\rm im}(f_s)$ is $r$ (the
rank of $E$).}
\end{definition}

It is easy to see that $H^0(\Sigma,\, E^{\oplus n})_0$ is a Zariski
open subset of $H^0(\Sigma,\, E^{\oplus n})$ (but it can be empty).
Note that it corresponds to the subset of holomorphic sections defining
stable $n$-pairs $(E,s)$, as described in Section~\ref{stabHitKob}:
if $s\, \in\, H^0(\Sigma,\, E^{\oplus n})_0$, then the quotient $E/{\rm im}(f_s)$
is either zero, or it is a torsion sheaf supported
at finitely many points.

Take any $s\,\in\, H^0(\Sigma,\, E^{\oplus n})_0$. Let
$$
{\mathcal K}\, :=\, {\rm ker}(f_s)\, \subset\,
{\mathcal O}^{\oplus n}_\Sigma
$$
be the kernel of the homomorphism $f_s$ in \eqref{e2}. Consider the dual homomorphism
$$
({\mathcal O}^{\oplus n}_\Sigma)^*\,\cong \,
{\mathcal O}^{\oplus n}_\Sigma\, \longrightarrow\, {\mathcal K}^*
$$
to the inclusion map ${\mathcal K}\,\hookrightarrow\, {\mathcal 
O}^{\oplus n}_\Sigma$. So ${\mathcal K}^*$ is a quotient bundle of
${\mathcal O}^{\oplus n}_\Sigma$.

We have a short exact sequence of coherent sheaves on $\Sigma$
\begin{equation}\label{e3a}
0\,\longrightarrow\, E^*\,\stackrel{f^*_s}{\longrightarrow}\,
{\mathcal O}^{\oplus n}_\Sigma\,\longrightarrow\,
{\mathcal K}^*\oplus {\mathcal T}\, =: {\mathcal Q}
\,\longrightarrow\, 0\, ,
\end{equation}
where $\mathcal T$ is either a torsion sheaf supported on finitely
many points of $\Sigma$, or ${\mathcal T}\,=\, 0$; in fact, $\mathcal T$
is isomorphic to the quotient sheaf $E/{\rm im}(f_s)$ (but there
is no canonical isomorphism).

Since $E^*$ is a subsheaf of a trivial vector bundle,
it follows that the degree of $E^*$ is never positive; hence we will require throughout that
\begin{equation}\label{de}
d:=\text{deg}(E)\, =\, -\text{deg}(E^*) \, \geq\, 0\, .
\end{equation}

We now introduce an ample line bundle ${\mathcal L} \longrightarrow \Sigma$ over
the surface where the vortices live. (For the purposes of the present section, this line bundle plays an 
auxiliary role, and its choice does not affect any of the results.)
Since $\ell := \text{deg}({\mathcal L})$ is necessarily positive, there is an 
integer
\begin{equation}\label{dE}
\delta_E\, \in\, {\mathbb N}
\end{equation}
such that, for all $\delta\, \geq\, \delta_E$,
\begin{equation}\label{e4}
H^1(\Sigma,\, E^*\otimes {\mathcal L}^{\otimes \delta})\, =\, 0\, ,
\end{equation}
and the natural evaluation homomorphism
\begin{equation}\label{e4a}
H^0(\Sigma,\, E^*\otimes{\mathcal L}^{\otimes
\delta})\otimes_{\mathbb C} {\mathcal O}_\Sigma\,\longrightarrow
\, E^*\otimes {\mathcal L}^{\otimes \delta}
\end{equation}
is surjective. The second condition means that the vector
bundle $E^*\otimes {\mathcal L}^{\otimes \delta}$
is generated by its global holomorphic sections.
We emphasize that at this stage $\delta_E$ depends on the 
holomorphic vector bundle $E$. The Riemann--Roch theorem yields
\begin{equation}\label{rr}
\dim H^0(\Sigma,\, E^*\otimes {\mathcal L}^{\otimes \delta})-
\dim H^1(\Sigma,\, E^*\otimes {\mathcal L}^{\otimes \delta})
\,=\, r \ell \delta  - d +r(1-g)\,
\end{equation}
and this determines the dimension of $H^0(\Sigma,\,E^*\otimes {\mathcal 
L}^{\otimes \delta})$ whenever
$\delta\ge \delta_E$, by \eqref{e4}.

Suppose that an integer  $\delta$ is fixed, satisfying $\delta\, \geq\, \delta_E$. Tensoring
\eqref{e3a} with ${\mathcal L}^{\otimes \delta}$, we obtain the
short exact sequence of coherent sheaves
\begin{equation}\label{e3}
0\,\longrightarrow\, E^*\otimes 
{\mathcal L}^{\otimes \delta}\,\longrightarrow\,
({\mathcal L}^{\otimes \delta})^{\oplus n}\
\longrightarrow\, {\mathcal Q}\otimes {\mathcal L}^{\otimes \delta}
\,\longrightarrow\, 0\, .
\end{equation}
This will give rise to a long exact sequence of cohomology groups
\begin{equation}\label{e5}
0\,\longrightarrow\, H^0(\Sigma,\, E^*\otimes{\mathcal L}^{\otimes
\delta}) \,\longrightarrow\, H^0(\Sigma,\, ({\mathcal L}^{\otimes
\delta})^{\oplus n})\, \stackrel{Q}{\longrightarrow}\,H^0(\Sigma,\,{\mathcal Q}
\otimes {\mathcal L}^{\otimes\delta}) \,\longrightarrow\, 0\, ,
\end{equation}
where the right-exactness follows from \eqref{e4}.

Consider the quotient map $Q: H^0(\Sigma,\, ({\mathcal L}^{\otimes
\delta})^{\oplus n})\, \longrightarrow\,H^0(\Sigma,\,{\mathcal Q}
\otimes {\mathcal L}^{\otimes\delta})$
in \eqref{e5}. The subsheaf
$$
E^*\otimes
{\mathcal L}^{\otimes \delta}\, \subset\,
({\mathcal L}^{\otimes \delta})^{\oplus n}
$$
in \eqref{e3} can be reconstructed from $Q$, and from it the
morphism $f_s$ in (\ref{e2}), by a 
procedure that we will now describe.

Let
$$\widehat{\mathcal K}\, :=\, \text{ker} \, Q$$
be the kernel of the quotient map, and let
\begin{equation}\label{S}
{\mathcal S}\, \subset\, ({\mathcal L}^{\otimes
\delta})^{\oplus n}
\end{equation}
be the subsheaf generated by the sections lying in the
subspace $\widehat{\mathcal K}$. From the exactness of
the sequence \eqref{e5} we know that $\widehat{\mathcal K}$
coincides with the subspace
$$
H^0(\Sigma,\, E^*\otimes{\mathcal L}^{\otimes
\delta}) \,\hookrightarrow\, H^0(\Sigma,\, ({\mathcal L}^{\otimes
\delta})^{\oplus n})\, 
$$
determined by the section $s \in H^0(\Sigma,\,E^{\oplus n})_0$. Also, the
holomorphic vector bundle $E^*\otimes{\mathcal L}^{\otimes
\delta}$ is generated by its global sections (recall that the
homomorphism in \eqref{e4a} is surjective). Consequently, the
subsheaf $\mathcal S$ in \eqref{S} coincides with the subsheaf
$$
E^*\otimes
{\mathcal L}^{\otimes \delta}\, \subset\,
({\mathcal L}^{\otimes \delta})^{\oplus n}
$$
in \eqref{e3}. In other words, we have reconstructed the
subsheaf $E^*\otimes{\mathcal L}^{\otimes \delta}$
of $({\mathcal L}^{\otimes \delta})^{\oplus n}$ from the
quotient vector space $H^0(\Sigma,\, ({\mathcal 
L}^{\otimes\delta})^{\oplus n})/H^0(\Sigma,\,E^*\otimes  {\mathcal 
L}^{\otimes\delta})$, or equivalently from the linear map $Q$ in (\ref{e5}).

Let $E'\, \longrightarrow\, \Sigma$ be a holomorphic
vector bundle, and let $s\, \in\, H^0(\Sigma,\, (E')^{\oplus n})_0$
(see Definition \ref{def1}) be such that
$$
H^1(\Sigma,\, (E')^*\otimes {\mathcal L}^{\otimes \delta})\, =\, 0\, ,
$$
and also assume that $(E')^*\otimes {\mathcal L}^{\otimes \delta}$ is 
generated
by its global sections. Let ${\mathcal Q}'$ be the quotient
of ${\mathcal O}^{\oplus n}_\Sigma$ constructed from $E'$ just as
$\mathcal Q$ is constructed from $E$ (see \eqref{e3a}). Therefore,
${\mathcal Q}'\otimes {\mathcal L}^{\otimes \delta}$ is a quotient
of $({\mathcal L}^{\otimes \delta})^{\oplus n}$.
If the two quotients
$H^0(\Sigma,\,{\mathcal Q}\otimes {\mathcal L}^{\otimes\delta})$ and
$H^0(\Sigma,\,{\mathcal Q}'\otimes {\mathcal L}^{\otimes\delta})$ of
$H^0(\Sigma,\, ({\mathcal L}^{\otimes
\delta})^{\oplus n})$ coincide, then the subsheaf
$E^*\otimes
{\mathcal L}^{\otimes \delta}$ of $({\mathcal L}^{\otimes 
\delta})^{\oplus n}$ (see \eqref{e3}) coincides with the subsheaf
$(E')^*\otimes
{\mathcal L}^{\otimes \delta}$ constructed as in \eqref{e3}
using $E'$. Indeed, this follows from the above observation
that we can reconstruct the subsheaf $E^*\otimes
{\mathcal L}^{\otimes \delta}$ of
$({\mathcal L}^{\otimes \delta})^{\oplus n}$ from the quotient map 
$Q: H^0(\Sigma,\, ({\mathcal L}^{\otimes
\delta})^{\oplus n})\, \longrightarrow\,H^0(\Sigma,\,{\mathcal Q}
\otimes {\mathcal L}^{\otimes\delta})$
in \eqref{e5}.

We put down the observations above in the form of what we will call a ``reconstruction'' lemma:

\begin{lemma}\label{lem1}
The quotient $H^0(\Sigma,\, ({\mathcal L}^{\otimes
\delta})^{\oplus n})\, \longrightarrow\,H^0(\Sigma,\,{\mathcal Q}
\otimes {\mathcal L}^{\otimes\delta})$
in \eqref{e5} uniquely determines the subsheaf
$$
E^*\otimes
{\mathcal L}^{\otimes \delta}\, \subset\,
({\mathcal L}^{\otimes \delta})^{\oplus n}
$$
in \eqref{e3}.
\end{lemma}

\begin{remark}\label{re-c}
{\rm Consider the subsheaf $E^*\otimes
{\mathcal L}^{\otimes \delta}\, \subset\,
({\mathcal L}^{\otimes \delta})^{\oplus n}$ in
Lemma \ref{lem1}. Its dual $E\otimes
({\mathcal L}^*)^{\otimes \delta}$ is a quotient
of $(({\mathcal L}^{\otimes \delta})^{\oplus n})^*
\,=\, (({\mathcal L}^*)^{\otimes \delta})^{\oplus n}$.
Tensoring this quotient homomorphism}
$$
(({\mathcal L}^*)^{\otimes \delta})^{\oplus n}
\, \longrightarrow\, E\otimes
({\mathcal L}^*)^{\otimes \delta}
$$
{\rm with the identity homomorphism of
${\mathcal L}^{\otimes \delta}$, we get back
the homomorphism}
$$
f_s\, :\, ({\mathcal O}_\Sigma)^{\oplus n}\,\longrightarrow\, E
$$
{\rm used to construct the quotient in \eqref{e5}. So the quotient effectively 
determines the $n$-pair.}
\end{remark}

\begin{proposition}\label{prop1}
Fix a positive integer $r$, a nonnegative integer $d$
and an integer $n\, \geq\, r$. Given an ample line bundle ${\mathcal L} 
\longrightarrow \Sigma$, there is an integer
$\delta_{n,r,d}$ such that for any $n$-pair
$(E\, ,s)$ with
${\rm rk}(E)=r$, ${\rm deg}(E)=d$ and
$$
s\, \in\, H^0(\Sigma,\, E^{\oplus n})_0
$$
(see Definition \ref{def1}), and any integer $\delta\, \geq\,
\delta_{n,r,d}$,
\begin{itemize}
\item the homomorphism in \eqref{e4a} is surjective,
and
\item \eqref{e4} holds.
\end{itemize}
\end{proposition}

\begin{proof}
The strategy of the
proof is to first show, using the idea of Quot scheme, that all
such pairs of the given numerical type form a bounded family; then
the proof is completed
using upper semicontinuity for dimension of cohomology.

Take a pair $(E\, ,s)$, where $E\, \longrightarrow\, \Sigma$ is a
holomorphic vector bundle of rank $r$ and degree $d$, and
$$
s\, \in\, H^0(\Sigma,\, E^{\oplus n})_0\, .
$$
The vector bundle $E^*$ is a subsheaf
of ${\mathcal O}^{\oplus n}_\Sigma$ of rank $r$ and degree
$-d$ (see \eqref{e3a}).
Therefore, all possible pairs
$(E^*\, ,f^*_s)$ (see \eqref{e3a}) are parametrized by a projective
scheme $\mathbb T$ over $\mathbb C$
of finite type \cite[p. 40, Theorem 2.2.4]{HL}
(set $S$ in \cite[Theorem 2.2.4]{HL} to be a point). Now from upper
semicontinuity of dimension of $H^1$ we conclude that there
is an integer $k_0$, that depends only on $n$, $r$ and $d$, such
that for all $(E\, ,s)$ of the above type and all $\delta\, \geq\, k_0$, 
$$
H^1(\Sigma,\, E^*\otimes {\mathcal L}^{\otimes \delta})\, =\, 0\, .
$$

Take any point $x\, \in\, \Sigma$.
Consider the short exact sequence of sheaves
$$0\, \longrightarrow\,E^*\otimes{\mathcal L}^{\otimes\delta}\otimes
{\mathcal O}_\Sigma(-x)\, \longrightarrow\,E^*\otimes{\mathcal 
L}^{\otimes\delta}\, \longrightarrow\,
(E^*\otimes {\mathcal L}^{\otimes \delta})_x
\, \longrightarrow\, 0\, .
$$
Let
\begin{equation}\label{le}
H^0(\Sigma,\, E^*\otimes{\mathcal L}^{\otimes
\delta})\, \longrightarrow\, (E^*\otimes {\mathcal L}^{\otimes 
\delta})_x \, \longrightarrow\,
H^1(\Sigma,\, E^*\otimes{\mathcal L}^{\otimes\delta}\otimes
{\mathcal O}_\Sigma(-x))
\end{equation}
be the corresponding long exact sequence in cohomology. From
\eqref{le} we conclude that if
\begin{equation}\label{l2}
H^1(\Sigma,\, E^*\otimes{\mathcal L}^{\otimes\delta}\otimes
{\mathcal O}_\Sigma(-x))\,=\, 0\, ,
\end{equation}
then the homomorphism $H^0(\Sigma,\, E^*\otimes{\mathcal L}^{\otimes
\delta})\, \longrightarrow\, (E^*\otimes {\mathcal L}^{\otimes
\delta})_x$ in \eqref{le} is surjective.
Therefore, given $(E^*\, ,f^*_s)$, if \eqref{l2} holds for all
$x\, \in\, \Sigma$, then the homomorphism in
\eqref{e4a} is surjective.

Now all possible pairs
$(E^*\, ,f^*_s)$ (see \eqref{e3a}) are parametrized by a projective
scheme over $\mathbb C$ (see above). From upper
semicontinuity of dimension of $H^1$, we conclude again that there
is an integer $k_1$ such that, for all $(E\, ,s)$ of the type
in the statement of the proposition, and
all $\delta\, \geq\, k_1$, the homomorphism in 
\eqref{e4a} is surjective.

Consequently, the integer
\begin{equation}\label{dc}
\delta_{n,r,d}\, :=\, \text{max}\{k_0\, ,k_1\}\, ,
\end{equation}
which depends only on $r$, $d$ and $n$,
has the property that for all $\delta\, \geq\,
\delta_{n,r,d}$, and for any pair $(E\, ,s)$ of the
of the type in the statement of the proposition,
the homomorphism in \eqref{e4a} is surjective, and \eqref{e4}
holds. This completes the proof of the proposition.
\end{proof}

Note that in Proposition \ref{prop1} we assume the degree $d$ to be nonnegative
because of the inequality in \eqref{de}.

\section{Embedding in a Grassmannian}\label{embed}

As in Section~\ref{sec1}, fix a positive integer $r$, a nonnegative
integer $d$ and an integer $n\, \geq\, r$, specifying the topology of 
$E \longrightarrow \Sigma$ 
and the number of copies of $E$ in a direct sum. For a given ample line 
bundle ${\mathcal L} \longrightarrow \Sigma$
of degree $\ell$, fix also an integer
$\delta\, \geq\, \delta_{n,r,d}$, where $\delta_{n,r,d}$ is
as in Proposition \ref{prop1}, cf.~\eqref{dc}. Notice that we can always set $\delta$
to be the minimal $\delta_{n,r,d}$ ensuring that both (\ref{e4}) is 
surjective
and the vanishing in (\ref{e4a}) holds, and in fact we will be doing so 
by default. A consequence of our previous discussion is that
\begin{equation}\label{inequal}
\ell \delta \,\ge\, \frac{d}{r}+g-1\, ;
\end{equation}
this follows from equations (\ref{rr}) and (\ref{e4}).

At this point, we shall introduce metric structures on 
the basic objects that we have been considering in the previous section.
We equip $\Sigma$ with a K\"ahler metric $\omega_\Sigma$, and the ample line 
bundle ${\mathcal L} \longrightarrow \Sigma$ in Section~\ref{sec1} with a Hermitian
structure $h_{\mathcal L}$. If the K\"ahler class $[ \omega_\Sigma] \in 
H^2(\Sigma,\,\mathbb{R})$ is integral, which amounts to
\[
\int_\Sigma \omega_\Sigma \,\in\, \mathbb{Z}\, ,
\]
it would be natural to require 
$({\mathcal L}, h_{\mathcal L})$, together with its Chern connection $\nabla_{\mathcal L}$, to be a prequantization of the K\"ahler structure on
$\Sigma$, in the sense that its curvature is proportional to the K\"ahler form as
\begin{equation} \label{prequant}
{F}_{\nabla_{\mathcal L}} \,=\, 2 \pi \sqrt{-1} \cdot \omega_{\Sigma}\, 
;
\end{equation}
but for now we need not impose this condition. 
Consider the vector space
$$H^0(\Sigma,\, ({\mathcal L}^{\otimes \delta})^{\oplus n})
\,\cong \, H^0(\Sigma,\, {\mathcal L}^{\otimes \delta})^{\oplus n}\, .$$
The Hermitian structure $h_{\mathcal L}$ on
$\mathcal L$ together with the K\"ahler form
$\omega_\Sigma$ on $\Sigma$ define an $L^2$ inner product on
$H^0(\Sigma,\, ({\mathcal L}^{\otimes \delta})^{\oplus n})$.

Let
\begin{equation}\label{gr}
\text{Gr}\, :=\, \text{Gr}(H^0(\Sigma,\, ({\mathcal L}^{\otimes
\delta})^{\oplus n})\, , r(\ell \delta -g +1)-d)
\end{equation}
be the Grassmannian of subspaces of
$H^0(\Sigma,\, ({\mathcal L}^{\otimes \delta})^{\oplus n})$
of dimension $r(\ell\delta-g+1)-d$ (see \eqref{rr} and \eqref{e5}).
The inner product on $H^0(\Sigma,\, ({\mathcal L}^{\otimes
\delta})^{\oplus n})$ defines a Fubini--Study K\"ahler
form on $\text{Gr}$. Indeed, for any subspace
$$
H^0(\Sigma,\, ({\mathcal L}^{\otimes \delta})^{\oplus n})
\, \supset\, V\, \in \, \text{Gr}\, ,
$$
we have
$$
T_V \text{Gr}\,=\,  V^*\otimes (H^0(\Sigma,\,
({\mathcal L}^{\otimes\delta})^{\oplus n})/V)\, , 
$$
where $T_V \text{Gr}$ is the holomorphic tangent space
at the point $V$ of $\text{Gr}$. The $L^2$ inner product we have
on $H^0(\Sigma,\, ({\mathcal L}^{\otimes\delta})^{\oplus n})$
defined above induces inner products on both $V$ and
$H^0(\Sigma,\, ({\mathcal L}^{\otimes\delta})^{\oplus n})/V$.
Therefore, we get an inner product on
$T_V \text{Gr}$. It is easy to see that the Hermitian
structure on $\text{Gr}$ constructed in this way is
actually K\"ahler.

Another way of describing this K\"ahler structure is to
consider the Fubini--Study metric on the projective space of lines in
${\bigwedge}^{r(\ell\delta-g+1)-d} H^0(\Sigma,\,( {\mathcal L}^{\otimes \delta})^{\oplus n})$
\[
{\mathbb P}({\wedge}^{r(\ell\delta-g+1)-d} H^0(\Sigma,\, ({\mathcal L}^{\otimes
\delta})^{\oplus n}))
\]
induced by the inner product on $H^0(\Sigma,\, {\mathcal L}^{\otimes \delta})$.
The Pl\"ucker map~\cite{GH}
\begin{equation} \label{Plu}
P\,:\, {\rm Gr}\, \longrightarrow\, 
{\mathbb P}(\bigwedge\nolimits^{r(\ell\delta-g+1)-d} H^0(\Sigma,\, 
({\mathcal 
L}^{\otimes
\delta})^{\oplus n})),
\end{equation}
defined by
$$
{\rm Gr}\; \ni\; {\rm span}_{\mathbb{C}}\{ s_1, \ldots, s_{r(\ell\delta -g+1)-d} \}\, \longmapsto \,
s_1\wedge \cdots \wedge s_{r(\ell \delta-g+1)-d},
$$
embeds $\text{Gr}$ as a complex submanifold of the target. 
The above K\"ahler structure on $\text{Gr}$
coincides with the restriction of the Fubini--Study metric on
the projective space to the image of $P$.

Let
\begin{equation}\label{mo}
{\mathcal M}_\Sigma \,:=\, {\mathcal M}_\Sigma(n,r,d)
\end{equation}
be the moduli space of isomorphism classes of all $n$-pairs 
$(E\, ,s)$, on $\Sigma$ where the holomorphic bundle $E\, \longrightarrow\, \Sigma$ has
rank $r$ and degree $d$, and
$$
s\, \in\, H^0(\Sigma,\, E^{\oplus n})_0\, .
$$
We now claim that we have an embedding
\begin{equation}\label{vp}
\varphi\, :\, {\mathcal M}_\Sigma\, \longrightarrow\,\text{Gr}
\end{equation}
that sends any $(E\, ,s)$ to the subspace
$H^0(\Sigma,\, E^*\otimes{\mathcal L}^{\otimes
\delta}) \,\subset\, H^0(\Sigma,\, ({\mathcal L}^{\otimes
\delta})^{\oplus n})$ in \eqref{e5}, where
$\text{Gr}$ is defined in \eqref{gr}. Note that
\eqref{e4} and \eqref{rr} together imply that
$H^0(\Sigma,\, E^*\otimes{\mathcal L}^{\otimes\delta})$ has dimension 
$r(\ell\delta - g + 1)-d$, and this means that $\varphi$ is well 
defined. The map $\varphi$ is also injective from Lemma~\ref{lem1} and 
Remark~\ref{re-c}. In this way, the moduli space ${\mathcal M}_\Sigma$
can be regarded as a closed subvariety of the Grassmannian ${\rm Gr}$ in (\ref{gr}). 

One advantage of our description of the moduli space ${\mathcal M}_\Sigma$ is that 
one can address its smoothness in a straightforward way.
Take any point $\underline{z}\, :=\,
(E\, ,s)\,\in\, {\mathcal M}_\Sigma$ of the moduli space. Let
$$
0\,\longrightarrow\, E^*\,\stackrel{f^*_s}{\longrightarrow}\,
{\mathcal O}^{\oplus n}_\Sigma\,\longrightarrow\,
{\mathcal K}^*\oplus {\mathcal T}\, =: {\mathcal Q}
\,\longrightarrow\, 0
$$
be the short exact sequence constructed in \eqref{e3a} from the $n$-pair 
$(E\, ,s)$. The tangent space to ${\mathcal M}_\Sigma$ at the point 
$\underline{z}\, :=\, (E\, ,s)$ has the following description:
\begin{equation}\label{tp}
T_{\underline{z}}{\mathcal M}_\Sigma\, =\, H^0(\Sigma,\,
{\mathcal H}om(E^*\, ,{\mathcal Q}))\, =\,
H^0(\Sigma,\, E\otimes {\mathcal Q})\, .
\end{equation}
The obstruction to smoothness of ${\mathcal M}_\Sigma$ at
$\underline{z}$ is given by
$$
\text{Ext}^1_{{\mathcal O}_\Sigma}(E^*\, , {\mathcal Q})\, ,
$$
where $\text{Ext}^1_{{\mathcal O}_\Sigma}$ is the global Ext.
Since $E^*$ is a vector bundle,
\begin{equation}\label{t2}
\text{Ext}^1_{{\mathcal O}_\Sigma}(E^*\, , {\mathcal Q})\,=\,
H^1(\Sigma,\, E^*\otimes {\mathcal Q}^*)\, . 
\end{equation}

In the local case, where $n = r:= \text{rank}(E)$, the quotients ${\mathcal Q}$
in \eqref{e3a} are torsion sheaves supported on finitely many
points of $\Sigma$. In that case, $E^*\otimes {\mathcal Q}^*$ is a 
torsion sheaf, and hence
$$
H^1(\Sigma,\, E\otimes {\mathcal Q})\, =\, 0\, .
$$
Therefore, from \eqref{t2} we conclude $\text{Ext}^1_{{\mathcal 
O}_\Sigma}(E^*\, , {\mathcal Q})\,=\,0$ if $n = r$, implying that
the variety ${\mathcal M}_\Sigma$ is smooth if $n = r$.

Since the map $\varphi$ in (\ref{vp}) embeds ${\mathcal M}_\Sigma$ 
in $\text{Gr}$ as a complex
submanifold, we can obtain K\"ahler structures on the 
moduli space ${\mathcal M}_\Sigma$ by
restricting a K\"ahler forms on the Grassmann manifold $\text{Gr}$ to it.
In the following, we shall denote by
$\omega_{\rm Gr}$ the K\"ahler form on the moduli space  
${\mathcal M}_\Sigma$ obtained by pulling back the Fubini--Study 
2-form on $\text{Gr}$ described above, using the embedding $\varphi$.
In the next section, we will see when it will be possible to make 
$\omega_{\rm Gr}$ cohomologous
to the usual K\"ahler structure $\omega_{L^2}$ on the moduli space of vortices,
in the abelian case where the K\"ahler class $[\omega_{L^2}]$ is known.

Although the K\"ahler form $\omega_{L^2}$ depends on both the metric on $\Sigma$ and a Hermitian metric on the vector bundle $E \longrightarrow \Sigma$, there is a natural splitting $\omega_{L^2}=\omega_1 + \omega_2$, where $\omega_1$ is a closed $(1,1)$-form depending
only on $\omega_\Sigma$
(see~\cite{ManNasVVMS} for the abelian case). A natural question to ask 
is 
how $\omega_1$ is related to $\omega_{\rm Gr}$ when the prequantization
condition (\ref{prequant}) is imposed. This is one issue that we plan to
address in future work.

\section{Abelian local vortices: $n=r=1$}\label{abelianlocal}

It is natural to ask whether the K\"ahler form $\omega_{\rm Gr}$ on the 
vortex moduli space ${\mathcal M}_\Sigma$
induced from the embedding $\varphi$ into the Grassmannian manifold, as described in
Section~\ref{embed}, is related to the $L^2$ K\"ahler structure 
inherited naturally from the gauge theory, 
which is of interest to physicists. 
We shall address this issue in the
present section, but our discussion will be restricted to the case of abelian
local vortices, where $n=1$, $r=1$. So throughout this section we will be assuming that
\[
{\mathcal M}_\Sigma \,= \, {\mathcal M}_\Sigma (1,1,d)\, .
\]

\subsection{Some background on the geometry of the abelian local case} 
\label{bgabellocal}

Let us briefly recall how the K\"ahler structures $\omega_{L^2}$ on ${\mathcal M}_\Sigma$ arise
in the abelian local case. There are many alternative descriptions of 
the
$L^2$ metrics of vortices, but here we will concentrate on a particularly insightful one given by Garc\'\i a-Prada in~\cite{GaDEP}, which uses infinite-dimensional symplectic geometry. The space of fields appearing as
variables in the vortex equations (\ref{vor1})--(\ref{vor2}) is ${\mathcal A}\times {\mathcal C}$, where
${\mathcal A}$ is the space of unitary connections on the line bundle $E 
\longrightarrow \Sigma$ and
${\mathcal C}=C^{\infty}(\Sigma,E)$ is the vector space of smooth 
sections of this bundle. Any two connections differ by a global real 
1-form on $\Sigma$ with values on the Lie algebra $\mathfrak{u}(1) \cong 
{\sqrt{-1}} \cdot \mathbb{R}$, so ${\mathcal A}$ is an affine space modelled 
on the vector space
$\Omega^1(\Sigma, \mathbb{R})$. Thus in fact ${\mathcal 
A}\times {\mathcal C}$ is a complex manifold
with complex structure induced from the one on $\Sigma$:
\begin{equation}\label{cxstr}
(\dot A, \dot \phi) \,\longmapsto\, (\ast A\, , \sqrt{-1} \, \dot{\phi})\, .
\end{equation}
The component of this map in the first factor is the Hodge star operator on $\Sigma$ acting on 1-forms, which squares to
$-{\rm id}_{\Omega^1(\Sigma)}$, whereas the component in the second 
factor is the complex structure on the fibers of the bundle $E 
\longrightarrow \Sigma$.

There is an action of the gauge group ${\rm Aut}_\Sigma (E)\cong C^{\infty}(\Sigma, {\rm U}(1))$
on fields $(A,\phi) \in {\mathcal A} \times {\mathcal C}$, namely
\begin{equation}\label{gaugetr}
(A,\phi) \,\longmapsto \,(A-u^{-1}{\rm d}u, u \phi)\, ,
\end{equation}
where $u \in {\rm Aut}_\Sigma (E)$. This action turns out to be  Hamiltonian with respect to a natural product symplectic form,
\begin{equation}\label{prodomega}
\omega_{\mathcal A} + \omega_{\mathcal C}\, ,
\end{equation}
defined on the space of fields. The factor denoted by  $\omega_{\mathcal 
A}$  in (\ref{prodomega}) is the Atiyah--Bott structure~\cite{AtiBot} on the space of connections $\mathcal A$, while $\omega_{\mathcal C}$ is the natural
symplectic  structure (of constant coefficients, hence closed) on $\mathcal C$ produced out of the 
K\"ahler structure on $\Sigma$ and the Hermitian metric on $E$. The latter is usually simply called the
$L^2$ structure (on $\mathcal C$), since it is associated to the metric
\[
|| \dot\phi ||^2_{L^2} =\int_\Sigma (\dot\phi, \dot\phi)_{h_{E}} \omega_\Sigma
\]
defined for all sections $\dot\phi \in C^{\infty}(\Sigma, E)\cong T_{\phi}{\mathcal C}$, for any $\phi \in {\mathcal C}$.
The complex structure (\ref{cxstr}) on the space of fields ${\mathcal A}\times{\mathcal C}$ preserves (\ref{prodomega}),
so one can regard this space as a K\"ahler manifold. 

The first vortex equation (\ref{vor1}) is invariant under the complex structure (\ref{cxstr}), so the infinite-dimensional submanifold ${\mathcal N}$ of solutions to this
equation (pairs $(A,\phi)$ where $\phi$ is a holomorphic section for the 
holomorphic structure on $E\longrightarrow \Sigma$ associated to the 
connection $A$, cf.~\cite{DonKro}) has an induced K\"ahler structure, which is
again preserved by the ${\rm Aut}_\Sigma (E)$-action (\ref{gaugetr}). It turns out that the left-hand side of the second vortex equation (\ref{vor2}) is a  moment map for this action. So the moduli space of solutions of both (\ref{vor1}) and (\ref{vor2}), where the action of the group of gauge transformations is quotiented out, can be understood as
the infinite-dimensional Meyer--Marsden--Weinstein quotient
\begin{equation}\label{MMWquotient}
{\mathcal M}_\Sigma = {\mathcal N}/\!\!/{\rm Aut}_\Sigma(E)\, .
\end{equation}
This receives  a symplectic structure,  denoted as $\omega_{L^2}$, and which is usually
referred to as the $L^2$ structure on the moduli space of vortices ${\mathcal M}_\Sigma$. In fact, this argument is formal, since we
are dealing with an infinite-dimensional quotient, but the intuitive picture just given is confirmed by the analysis
carried out in~\cite{BD1,BD2,GaDEP}, which is itself quite insightful. The K\"ahler form $\omega_{L^2}$ satisfies the properties 
\[
{p^*}\omega_{L^2} =i^*(\omega_{\mathcal A} + \omega_{\mathcal C}),
\qquad i:{\mathcal N} \hookrightarrow {\mathcal A}\times {\mathcal C}\, 
,
\]
where $p$ denotes the projection from ${\mathcal N}$ to the space of 
${\rm Aut}_\Sigma (E)$-orbits.

Under the stability condition (\ref{stable}),
Bradlow~\cite{Br} and Garc\'\i a-Prada~\cite{GaDEP} showed that the quotient
${\mathcal M}_\Sigma$ in (\ref{MMWquotient}) can be identified with 
the $d$-th symmetric power of
$\Sigma$,
\begin{equation}\label{sym}
{\mathcal M}_\Sigma \cong {\rm Sym}^d(\Sigma):=\Sigma^d/\mathfrak{S}_d
\end{equation}
as a complex manifold. 
This space parametrizes effective divisors of degree $d$, interpreted as portrayals of vortex
locations on $\Sigma$. But the symplectic structure $\omega_{L^2}$ on ${\mathcal M}_\Sigma$ turns out to be much more difficult to describe explicitly.

When comparing $\omega_L^2$ with $\omega_{\rm Gr}$, the most
basic question to ask is whether the two are cohomologous (up to a scalar multiple, say) for any choice of the
data. The answer to this question is trivially affirmative if $g=0$, since then ${\rm Sym}^d(\Sigma) \cong \mathbb{P}^d$ and
$H^2(\mathbb{P}^d,\,\mathbb{Z})\,\cong\, \mathbb{Z}$, so the interesting setting for this question is $g\ge 1$.
Then the cohomology ring of ${\rm Sym}^d(\Sigma)$ is more complicated; the intersection
\begin{equation}\label{intersection}
H^{1,1} ({\rm Sym}^d(\Sigma),\mathbb{C}) \cap H^2 ({\rm Sym}^d(\Sigma),\, \mathbb{Z})\, ,
\end{equation}
where the K\"ahler classes of the moduli space are contained, turns out to be a rank two
lattice.
The K\"ahler class $[\omega_{L^2}]$ has been computed 
as~\cite{ManNasVVMS, BatL2M}
\begin{equation}\label{Kclass}
[\omega_{L^2}]=\left( \pi \tau  {\rm Vol}(\Sigma) - \frac{4 \pi^2 
d}{e^2}\right) \eta + \frac{ 2\pi^2}{e^2} \sigma\, ;
\end{equation}
a description of the generators $\eta$ and $\sigma$ of (\ref{intersection}) will be
given in Section~\ref{representability}.
It is  remarkable that this
formula involves so little detail on the geometrical data needed to set up the vortex equations
and to define the $L^2$ metric. 
In the following, we shall be interested in calculating the K\"ahler class
 $[\omega_{\rm Gr}]$ and relating it with  $[\omega_{L^2}]$.
The result (\ref{Kclass}) has been used to compute the symplectic volume of the moduli
space~\cite{ManNasVVMS} and the total scalar curvature~\cite{BatL2M} --- such quantities carry only cohomological information.

\subsection{Description of the embedding} \label{descremb}

In the abelian local case, we can describe the embedding (\ref{vp}) 
constructed
in Section~\ref{embed} more explicitly. More precisely, we will be interested in characterizing the composition $P\circ \varphi$, where $P$ is the Pl\"ucker embedding
(\ref{Plu}).
Given the result (\ref{sym}), the map we are interested in is
\begin{equation}\label{Pvi}
P\circ \varphi\, :\,{\rm Sym}^d(\Sigma) \,\longrightarrow\,
{\mathbb P}({\wedge}^{\ell \delta-g-d+1} H^0(\Sigma,\, {\mathcal L}^{\otimes
\delta}))
\end{equation}
where $P$ and $\varphi$ are constructed in \eqref{Plu} and \eqref{vp},
respectively. We shall give a description the holomorphic line bundle on
${\rm Sym}^d(\Sigma)$ associated to this projective embedding.

Let $p_1$ (respectively, $p_2$) be the projection of
${\rm Sym}^d(\Sigma)\times \Sigma$ to ${\rm Sym}^d(\Sigma)$
(respectively, $\Sigma$). Let also
$$
\Delta_0\, \subset\, {\rm Sym}^d(\Sigma)\times \Sigma
$$
be the tautological divisor consisting of all points
$(z\, ,x)\, \in\, {\rm Sym}^d(\Sigma)\times \Sigma$ such that
$x\, \in\, z$.

Consider the line bundle $p^*_2
{\mathcal L}^{\otimes \delta}$ on ${\rm Sym}^d(\Sigma)\times \Sigma$,
and the torsion sheaf defined by
$$
{\mathcal B}\,:=\, p^*_2 {\mathcal L}^{\otimes \delta}/(p^*_2
{\mathcal L}^{\otimes \delta} \otimes {\mathcal O}_{{\rm
Sym}^d(\Sigma)\times \Sigma}(-\Delta_0))
\, \longrightarrow\, {\rm Sym}^d(\Sigma)\times \Sigma.
$$
The support of ${\mathcal B}$ is
$\Delta_0$, which is finite over ${\rm Sym}^d(\Sigma)$ of
degree $d$. Hence the direct image 
$$
p_{1*}{\mathcal B}\, \longrightarrow\, {\rm Sym}^d(\Sigma)
$$
is a vector bundle of rank $d$. So
$\bigwedge^d p_{1*}{\mathcal B}$ is a line bundle over
${\rm Sym}^d(\Sigma)$.

We have a canonical isomorphism of line bundles over
${\rm Sym}^d(\Sigma)$
\begin{equation}\label{c1}
\bigwedge\nolimits^d p_{1*}{\mathcal B}\,=\,
(P\circ\varphi)^*{\mathcal O}_{{\mathbb P}({\wedge}^{\ell-g-d+1}
H^0(\Sigma,\, {\mathcal L}^{\otimes\delta}))}(1)\, ,
\end{equation}
where $P\circ\varphi$ is the map in \eqref{Pvi}, and
$$
{\mathcal O}_{{\mathbb P}({\wedge}^{\ell-g-d+1}
H^0(\Sigma,\, {\mathcal L}^{\otimes\delta}))}(1)\,
\longrightarrow\,{\mathbb P}({\wedge}^{\ell-g-d+1}
H^0(\Sigma,\, {\mathcal L}^{\otimes\delta}))
$$
is the tautological line bundle. This means that the embedding \eqref{Pvi} is
associated to the complete linear system corresponding to the holomorphic
line bundle $\bigwedge^d p_{1*}{\mathcal B} \longrightarrow {\rm Sym}^d(\Sigma)$.

\subsection{Representability of the $L^2$ K\"ahler structure} 
\label{representability}

Our main goal in this section is to prove the following representability result:

\begin{theorem} \label{cohomL2}
Consider the embedding (\ref{Pvi}), constructed from an ample line 
bundle ${\mathcal L} \longrightarrow \Sigma$
of degree $\ell$ and an integer $\delta>\delta_{1,1,d}$, where $\delta_{1,1,d}$  is as in 
Proposition~\ref{prop1} and $d>1$. Then the Fubini--Study metric on 
${\rm Sym}^d(\Sigma)$ (obtained by pulling back the usual 
Fubini--Study metric using this map) is cohomologous to a multiple of the 
$L^2$-metric of vortices on the line bundle $E \longrightarrow \Sigma$ 
exactly when
\begin{equation}\label{quantiz}
q:=\frac{\tau e^2 }{4 \pi} {\rm Vol}(\Sigma)  \; \in \; \mathbb{N}
\end{equation}
and the integers $\ell, \delta$ are chosen such that
\begin{equation}\label{elldelta}
\ell \delta =q + g - 1\, .
\end{equation}
\end{theorem}

This result means that, at least in the abelian local case, the K\"ahler 
structure $\omega_{\rm Gr}$ on
${\mathcal M}_\Sigma$ discussed 
in Section~\ref{embed} provides a realization of the K\"ahler class
of the $L^2$ geometry of vortices if (\ref{quantiz}) and (\ref{elldelta}) hold. 
Note that the condition (\ref{quantiz}) is rather natural from the point of view of geometric
quantization, as it implies that the symplectic structure
$\frac{e^2}{2 \pi^2}\, \omega_{L^2}$
is (pre)quantizable in the sense of Weil:
\begin{equation}\label{weil}
\left[ \frac{e^2}{2 \pi^2}\, \omega_{L^2} \right] \in H^2(\mathcal{M}_\Sigma,\,\mathbb{Z})
 \end{equation} 
  (From (\ref{Kclass}), it follows that the Weil quantization
condition (\ref{weil}) is equivalent to $q \in \frac{1}{2}\mathbb{N}$.)  
It would be very striking if the full $L^2$ geometry were to be
described by a Fubini--Study structure, but we will not attempt to address this question here. Even in
the case $g=0$, for which the
representability of $[\omega_{L^2}]$ in the sense we are using is trivial, this question has not yet been settled rigorously.

To set the stage for the proof of Theorem~\ref{cohomL2}, we introduce the following curves on ${\rm Sym}^d (\Sigma)$, regarded as the space of
degree $d$ effective divisors on $\Sigma$:
\begin{eqnarray*}
\Sigma_\emptyset &:=& \{ d x \; | \;  x \in \Sigma  \}, \\
\Sigma_p &:=& \{ p + (d-1) x \; | \;  x \in \Sigma  \}, \qquad \text{for} \quad p \in \Sigma. 
\end{eqnarray*}
We shall denote their homology classes by
\[
\Sigma_0=[\Sigma_{\emptyset}] \qquad \text{and} \qquad  \Sigma_1=[\Sigma_p], 
\]
respectively. (Clearly, the homology class represented by $\Sigma_p$ is independent of $p\in \Sigma$ because we are assuming that $\Sigma$ is connected.) Let us also set
\begin{eqnarray}
d_0 &:=& {\rm deg} (P\circ \varphi (\Sigma_\emptyset)),\label{d0} \\
d_1 &:=& {\rm deg} (P\circ \varphi (\Sigma_p)) \label{d1}
\end{eqnarray}
to be the degrees of the images of the curves above by the map (\ref{Pvi}), whose target is a 
complex projective space of dimension
\[
N\,=\,\binom{ \ell \delta - g + 1}{d} -1\, .
\]

We first claim that the integers $d,g,d_0$ and $ d_1$ determine the cohomology class
\[
[(P\circ \varphi)^* \omega_{\rm FS}]
\,\in\, H^{1,1} ({\rm Sym}^d(\Sigma)) \cap H^2 ({\rm Sym}^d(\Sigma),\, \mathbb{Z})\, .
\]

To describe this, we start by recalling the basic result~\cite{McD}
\[
H^{k}({\rm Sym}^d (\Sigma),\,\mathbb{Z})\,\cong\, 
H^{k}(\Sigma^d,\,\mathbb{Z})^{\mathfrak{S}_d}, \qquad k\,\in\,\mathbb{N}\, 
.
\] 
The intersection of  cohomology groups in (\ref{intersection}) is generated over $\mathbb{Z}$ by the two cohomology classes of degree two~\cite{ManNasVVMS}
\begin{equation} \label{generators}
\eta \,=\,\sum_{i=1}^d \beta_i \qquad \text{and} \qquad 
\sigma\,=\,\sum_{j=1}^{g}\sigma_j\, .
\end{equation}
Here, the cohomology classes $\beta_i$ come from the fundamental class 
$ \beta \in H^2(\Sigma,\,\mathbb{Z})$; more precisely,
$\beta_i=\pi_i^*\beta$,
where $\pi_i: \Sigma^d \longrightarrow \Sigma$ denotes the projection to 
the $i$-th factor. Moreover, we denote
\begin{equation}\label{sigmaxi}
\sigma_j:=\xi_j  \xi_{j+g}, \qquad \text{where}\qquad \xi_j=\sum_{k=1}^{d} \alpha_{j,k},
 \qquad 1\le j \le 2g,
\end{equation} 
and 
the $\alpha_{j,k}$ are classes of degree one which come from the middle 
cohomology of $\Sigma$, namely
\[
\alpha_{j,k}=\pi_k^* \alpha_j\, .
\]
In this expression, the $\alpha_j$ denote elements in a standard basis of
$H^1(\Sigma,\,\mathbb{Z})$,
satisfying~\cite{Ful}
\begin{eqnarray*}
 &\alpha_i \alpha_j =0 & i \ne j \pm g, \\
 & \alpha_i \alpha_{i+g}=-\alpha_{i+g}\alpha_{i} = \beta & 1 \le i \le g.
\end{eqnarray*}
So we may write
\begin{equation} \label{pbFS}
(P \circ \varphi)^* [\omega_{\rm FS}] \,= \,C_\eta \eta + C_\sigma 
\sigma\, ,
\end{equation}
where $\eta$ and $\sigma$ are the generators in (\ref{generators}),
so our task is to obtain the coefficients $C_\eta, C_\sigma \in \mathbb{Z}$ as functions of $d, g, d_0$ and $d_1$.

\begin{lemma} \label{pairings}
The duality pairing on ${\rm Sym}^d(\Sigma)$ satisfies:
\[
\langle \eta, \Sigma_j \rangle = d-j \quad \text{and} \quad
\langle \sigma, \Sigma_j \rangle = (d-j)^2 g  \quad \text{for} \quad j\in \{ 0,1\}.
\]
\end{lemma}

\begin{proof}
The pairings for $j=0$ 
can be reduced to computations in the cohomology ring of ${\rm Sym}^d(\Sigma)$, which has
been given a presentation in~\cite[(6.3)]{McD}. In fact, the statement 
in reference~\cite{McD} is
not totally accurate --- we refer the reader to Section~2 of~\cite{BT} for the corrected result.
For our purposes, it will suffice to state that $H^*({\rm Sym}^d(\Sigma),\,\mathbb{Z})$ is generated
by the classes $\eta$ in (\ref{generators}) and $\xi_j$ in (\ref{sigmaxi}), $j=1,\ldots,2g$, which supercommute according to the parity of their degrees;
in particular, one has
\[
\eta\sigma_j=\sigma_j\eta, \qquad j=1,\ldots,g,
\]
where $\sigma_j$ were defined in (\ref{sigmaxi}), since $\eta$ and $\xi_j$ commute. The extra relations among the generators can be expressed as follows: given three disjoint subsets
\[
I_1,I_2, J \subset \{ 1,\ldots, g\}
\]
and a nonnegative integer $r$ satisfying~\cite[(2.3)]{BT}
\[
r \ge d-|I_1| - |I_2| - 2|J|+1,
\]
there is a nontrivial relation
\begin{equation}\label{nontrrel}
\eta^r \prod_{i_1 \in I_1} \xi_{i_1} \prod_{i_2 \in I_2} \xi_{i_2+g}\prod_{j\in J}(\eta-\sigma_j)=0.
\end{equation}

For $d\ge 1$, we have the relation
\begin{equation}\label{getadnew}
\eta^{d-1} \sigma = g \eta^d.
\end{equation}
This follows from summing the relations
\begin{equation} \label{sigjeta}
\eta^{d-1}\sigma_j =\eta^d,\qquad j=1,\ldots, g
\end{equation}
over $j$. Notice that (\ref{sigjeta}) can be obtained from 
(\ref{nontrrel}) by taking  $r=d-1$, $J=\{ j \}$ and $I_1=I_2=\emptyset$.

Another relation contained in~({\ref{nontrrel}}) is that, for $i\ne j$ and $d>1$,
\begin{equation}\label{sigisigj}
 \eta^{d-2} \sigma_i \sigma_j=\eta^{d-1} (\sigma_i + \sigma_j) -\eta^d;
\end{equation}
this one is obtained by setting $r=d-2, J=\{ i,j\}$ and 
$I_1=I_2=\emptyset$. Since $\sigma_j^2=0$ from the anticommutativity of 
the $\xi_j$'s (for each $j=1,\ldots, g$), we also have that
\begin{eqnarray}
\eta^{d-2} \sigma^2  &=& 2 \sum_{1\le i<j \le g}\eta^{d-2}  \sigma_i \sigma_j \nonumber\\
&=& 2 \sum_{1\le i<j \le g}\eta^{d-1} ( \sigma_i + \sigma_j ) - g (g-1) \eta^{d} \label{1ststep}\\
&=&g(g-1)\eta^{d}. \label{2ndstep}
\end{eqnarray}
Step (\ref{1ststep}) made use of (\ref{sigisigj}), whereas (\ref{2ndstep}) used (\ref{sigjeta}).

Another useful result by Macdonald~\cite[(15.4)]{McD} is that the Poincar\'e dual of the homology class $\Sigma_0$, for $d>1$, is given by
\begin{equation}\label{PDSig0}
{\rm PD}(\Sigma_0) = d(d+(g-1)(d-1))\eta^{d-1} - d(d-1)\eta^{d-2}\sigma.
\end{equation}
This can be applied to calculate
\begin{eqnarray}
\langle \eta, \Sigma_0 \rangle & 
= & d \int_{{\rm Sym}^d \Sigma}  (d+(d-1)(g-1))\eta^{d} - (d-1)\eta^{d-1}\sigma\nonumber \\
&=& d \int_{{\rm Sym}^d \Sigma}  (d+(d-1)(g-1) -(d-1)g)\eta^d \label{second}\\
&=&   d\int_{{\rm Sym}^d \Sigma}   \eta^d\nonumber\\
& =&d\label{fundamclass}.
\end{eqnarray}
The second step (\ref{second}) used (\ref{getadnew}), whereas the last 
step (\ref{fundamclass}) follows from the fact that $\eta^d$ is the 
fundamental class of ${\rm Sym}^d(\Sigma)$.

Using (\ref{PDSig0}) once again, we can write
\begin{eqnarray}
\langle \sigma, \Sigma_0 \rangle & 
= & d \int_{{\rm Sym}^d \Sigma}  (d+(d-1)(g-1))\sigma \eta^{d-1} - (d-1)\sigma \eta^{d-2}\sigma\nonumber \\
&=& d \int_{{\rm Sym}^d \Sigma}  ((d+(d-1)(g-1))g -(d-1)g(g-1))\eta^d\label{middle} \\
&=&  d^2 g, \nonumber
\end{eqnarray}
where (\ref{middle}) is a consequence of (\ref{sigjeta}) and (\ref{2ndstep}).

Now consider the map $\iota: \Sigma \longrightarrow \Sigma_p$ given by
\[
x \longmapsto p+ (d-1)x \quad  \in \quad {\rm Sym}^d(\Sigma)\, ,
\]
which is a biholomorphism.  We have 
\begin{equation}\label{iotaeta}
\iota^*\eta = (d-1) \beta
\end{equation}
and
\[
\iota^*\xi_j = (d-1) \alpha_j, \qquad j=1,\ldots, 2g,
\]
which in turn implies
\[
\iota^*({\xi_j}\xi_{j+g}) = (d-1)^2 \alpha_j\alpha_{j+g} = (d-1)^2 \beta, \qquad j=1,\ldots, g.
\]
It follows that
\begin{equation}\label{iotasigma}
\iota^*\sigma=(d-1)^2g \beta\, .
\end{equation}

So we can finally compute
\begin{eqnarray*}
\langle \eta, \Sigma_1\rangle& = &\int_{\Sigma_p}\eta \\
&=& \int_{\Sigma}\iota^*\eta\\
&=&(d-1)\int_\Sigma\beta\\
&=&d-1
\end{eqnarray*}
using (\ref{iotaeta}), and likewise, from (\ref{iotasigma}),
$$
\langle \sigma\, , \Sigma_1\rangle
\, =\, (d-1)^2 g\, .
$$
This completes the proof of the lemma.
\end{proof}

Since
\[
\langle (P \circ \varphi )^* [\omega_{\rm FS}] , \Sigma_j \rangle=d_j \qquad \text{for} \quad j=0,1,
\]
the constants $C_\eta$ and $C_\sigma$ in (\ref{pbFS}) can be determined by solving a linear system
whose coefficients are the four pairings in Lemma~\ref{pairings}. The solution is
\begin{eqnarray}
C_\eta & = & \frac{d^2 d_1 - (d-1)^2  d_0}{d(d-1)}, \label{Ceta}  \\ 
C_\sigma & = & \frac{(d-1)d_0 - d d_1} {d(d-1)g}\label{Csigma},
\end{eqnarray}
and this establishes our claim.

We want to compare the resulting K\"ahler class (\ref{pbFS}) with the  
K\"ahler 
class $[\omega_{L^2}]$ in (\ref{Kclass}) associated to the $L^2$ metric of vortices.
The next task is to determine the degrees $d_0$ and $d_1$ defined in (\ref{d0})--(\ref{d1}),
in terms of the basic topological data.

\begin{lemma} \label{d0d1}
$d_j=(d-j) (\ell \delta  + (d-j-1)(g-1)-j)$ for $j\in \{0,1\}$.
\end{lemma}

\begin{proof}
Let
\begin{equation}\label{j1}
\psi\, :\, \Sigma\, \longrightarrow\, {\rm Sym}^d (\Sigma)
\end{equation}
be the morphism defined by $x\, \longmapsto\, dx$. Note that
$d_0$ in \eqref{d0} is the degree of
$$
(P\circ \varphi\circ \psi)^*{\mathcal O}_{{\mathbb P}({\wedge}^{\ell 
\delta-g-d+1}H^0(\Sigma,\, {\mathcal L}^{\otimes
\delta}))}(1)\, ,
$$
where $P\circ \varphi$ is the morphism in \eqref{Pvi}.

Let $K_\Sigma$ be the holomorphic cotangent bundle of $\Sigma$.

Take any point $x\, \in\, \Sigma$. We have a natural filtration of
coherent sheaves
$$
{\mathcal L}^{\delta}\otimes {\mathcal O}_\Sigma(-dx)\,\subset\,
{\mathcal L}^{\delta}\otimes {\mathcal O}_\Sigma((1-d)x)\,\subset\,
\cdots\,\subset\,{\mathcal L}^{\delta}\otimes {\mathcal O}_\Sigma(-ix)
$$
$$
\subset \, {\mathcal L}^{\delta}\otimes {\mathcal O}_\Sigma((1-i)x)
\,\subset\, \cdots \,\subset\, {\mathcal L}^{\delta}\otimes {\mathcal 
O}_\Sigma(-x)\,\subset\, {\mathcal L}^{\delta}\, .
$$
For any $i\, \in\, [1\, ,d]$, the quotient ${\mathcal 
L}^{\delta}\otimes {\mathcal O}_\Sigma((1-i)x)/{\mathcal
L}^{\delta}\otimes {\mathcal O}_\Sigma(-ix)$ is the torsion
sheaf ${\mathcal L}^{\delta}_x\otimes (K^{\otimes (i-1)}_\Sigma)_x$
supported at $x$. Consequently, we have a canonical identification
$$
(P\circ \varphi\circ\psi)^*{\mathcal O}_{{\mathbb P}({\wedge}^{\ell
\delta-g-d+1} H^0(\Sigma,\,{\mathcal L}^{\otimes
\delta}))}(1)_x\,=\,
{\mathcal L}^{d\delta}_x\otimes (K^{\otimes d(d-1)/2}_\Sigma)_x
$$
(see \eqref{c1}), where $\psi$ is the map in \eqref{j1}.
Moving $x$ over $\Sigma$, this isomorphism produces
an isomorphism of line bundles
$$
(P\circ \varphi\circ\psi)^*{\mathcal O}_{{\mathbb P}({\wedge}^{\ell 
\delta-g-d+1} H^0(\Sigma,\, {\mathcal L}^{\otimes
\delta}))}(1)\,=\,
{\mathcal L}^{d\delta}_x\otimes (K^{\otimes d(d-1)/2}_\Sigma)
$$
Since $\text{deg}(K_\Sigma)\,=\, 2(g-1)$,
this immediately implies that $d_0\, =\, 
d(\ell\delta + (g-1)(d-1))$.

Fix a point $p\, \in\, \Sigma$. Let
\begin{equation}\label{j2}
\psi_1\, :\, \Sigma\, \longrightarrow\, {\rm Sym}^d (\Sigma)
\end{equation}
be the morphism defined by $x\, \longmapsto\, p+(d-1)x$. Note that
$d_1$ in \eqref{d1} coincides with
$$
\text{deg}((P\circ \varphi\circ \psi_1)^*{\mathcal O}_{{\mathbb 
P}({\wedge}^{\ell
\delta-g-d+1}H^0(\Sigma,\, {\mathcal L}^{\otimes
\delta}))}(1))\, ,
$$
where $P\circ \varphi$ is the morphism in \eqref{Pvi}.

For notational convenience, the line bundle
${\mathcal L}^{\delta}\otimes{\mathcal O}_\Sigma(-p)$ will
be denoted by $\zeta$.

As before, take any point $x\, \in\, \Sigma$.
We have a natural filtration of coherent sheaves
$$
\zeta\otimes {\mathcal O}_\Sigma((1-d)x)\,\subset\,
\zeta \otimes {\mathcal O}_\Sigma((2-d)x)\,\subset\,
\cdots\,\subset\, \zeta\otimes {\mathcal 
O}_\Sigma(-x)\,\subset\, \zeta\, .
$$
For any $i\, \in\, [1\, ,d-1]$, the quotient $\zeta\otimes {\mathcal 
O}_\Sigma((1-i)x)/{\mathcal
L}^{\delta}\otimes {\mathcal O}_\Sigma(-ix)$ is the torsion
sheaf $\zeta_x\otimes (K^{\otimes (i-1)}_\Sigma)_x$
supported at $x$. Consequently, we have a canonical identification
$$
(P\circ \varphi\circ\psi_1)^*{\mathcal O}_{{\mathbb P}({\wedge}^{\ell 
\delta-g-d+1} 
H^0(\Sigma,\, {\mathcal L}^{\otimes
\delta}))}(1)_x\,=\,
\zeta^{\otimes (d-1)}_x\otimes (K^{\otimes (d-1)(d-2)/2}_\Sigma)_x
\otimes {\mathcal L}^{\otimes \delta}_p\, ,
$$
(see \eqref{c1}), where $\psi_1$ is the map in \eqref{j2}.
Fixing an isomorphism of the line ${\mathcal L}^{\otimes \delta}_p$
with $\mathbb C$ (recall that $p$ is fixed), and moving $x$
over $\Sigma$, the above isomorphism gives an isomorphism
of line bundles
$$
(P\circ \varphi\circ\psi_1)^*{\mathcal O}_{{\mathbb P}({\wedge}^{\ell
\delta-g-d+1}
H^0(\Sigma,\, {\mathcal L}^{\otimes
\delta}))}(1)\,=\,
\zeta^{\otimes (d-1)}_x\otimes (K^{\otimes (d-1)(d-2)/2}_\Sigma)\, .
$$
Since $\text{deg}(\zeta)\,=\, \ell\delta-1$, this implies that
$d_1\,=\, (\ell\delta-1)(d-1) +(g-1)(d-1)(d-2)$.
\end{proof}

\begin{proof}[Proof of Theorem~\ref{cohomL2}] 
Using Lemma~\ref{d0d1} in (\ref{Ceta}) and  (\ref{Csigma}), we find
\begin{eqnarray*}
C_\eta & = & \ell \delta -d-g+1, \label{Cetanew}  \\ 
C_\sigma & = & 1\label{Csigmanew},
\end{eqnarray*}
which are integers as expected. 
Comparing with the coefficients of $\eta$ and $\sigma$ in (\ref{Kclass}),
the formula (\ref{elldelta}) for the quantity $q$ in Theorem~\ref{cohomL2} immediately follows. The quantization condition
(\ref{quantiz}) results from all the other terms in (\ref{elldelta}) being 
integers.
\end{proof}

Note that when the inequality (\ref{stable}) ensuring stability is saturated, in the situation
\begin{equation}\label{taucrit}
\tau \rightarrow \frac{4 \pi d}{e^2 {\rm Vol}(\Sigma)}
\end{equation}
which is called the limit of ``dissolved'' vortices by Manton and Rom\~ao in~\cite{ManRom}, the 
quantization condition~(\ref{quantiz}) is automatically satisfied with $q=d$. Then imposing the condition (\ref{elldelta}) implies that (\ref{inequal}) also becomes an equality, which unfortunately makes the Grassmannian
(\ref{gr}) collapse. The nontrivial situation closest to this collapse would be to consider
\[
\tau=\frac{4 \pi (d+1)}{e^2 {\rm Vol}(\Sigma)} \quad \Rightarrow \quad q=d+1,
\]
for which the Grassmannian (\ref{gr}) is a projective space; if the area of $\Sigma$ is taken to be large,
this value of $\tau$ will still be close to the critical value (\ref{taucrit}).
In this context (provided $\delta_{1,1,d}$ does not turn out to be too large), the geometry of the K\"ahler structure $\omega_{\rm Gr}$ we introduced in Section~\ref{embed}, assuming $\ell \delta=d+g$, should give 
an approximation of the $L^2$ geometry of
the moduli spaces, as an extension of the work by Baptista and
Manton~\cite{BatMan} in the case $g=0$.

\medskip
\noindent
{\bf Acknowledgements.}\, Both authors are grateful to the Harish-Chandra
Research Institute, Allahabad, for hosting a meeting where this work was
initiated. IB is also grateful to the Center for the
Topology and Quantization of Moduli Spaces (CTQM), Aarhus University, for hospitality.
The work of NMR was partially supported by CTQM, Aarhus University, and by the European
Commission in the framework of the Marie Curie project MTKD-CT-2006-042360; he would
like to thank Jo\~ao Baptista for several discussions.


\bibliographystyle{numsty}

\begin{thebibliography}{9999999}
                        \setlength{\itemsep}{0 pt}
                          \setlength{\parskip}{0pt}
                          \setlength{\parsep}{0pt}

\begin{small}

\bibitem[A]{At}
\textsc{M.F. Atiyah}: Complex analytic connections in fibre bundles.
\newblock \textsl{Trans.\ Amer.\ Math.\ Soc.} \textbf{85} (1957) 181 --207

\bibitem[AB]{AtiBot}
\textsc{M.F. Atiyah {\upshape and} R. Bott}: The Yang--Mills equations
over Riemann surfaces. 
\newblock \textsl{Phys.\  Trans.\ R.\ Soc,\ London} \textbf{308} (1982) 3495--3508

\bibitem[Ba1]{BatNAV}
\textsc{J.M. Baptista}: Non-abelian vortices on compact Riemann surfaces.
\newblock \textsl{Commun.\ Math.\ Phys.} \textbf{291} (2009) 799--812

\bibitem[Ba2]{BatL2M}
\textsc{J.M. Baptista}: On the $L^2$-metric of vortex moduli spaces.
\newblock \textsl{Nucl. Phys. B} \textbf{844} (2011) 308--333

\bibitem[BM]{BatMan}
\textsc{J.M. Baptista {\upshape and} N.S. Manton}: 
The dynamics of vortices on $S^2$ near the Bradlow limit.
\newblock \textsl{J.\ Math.\ Phys.} \textbf{44} (2003) 3495--3508

\bibitem[BDW2]{BDW}
\textsc{A. Bertram, G. Daskalopoulos {\upshape and} R. Wentworth}: 
Gromov invariants for holomorphic maps from Riemann surfaces to
Grassmannians.
\newblock \textsl{J. Amer. Math. Soc.} \textbf{9} (1996) 529--571

\bibitem[BT]{BT}
\textsc{A. Bertram {\upshape and} M. Thaddeus}: 
On the quantum cohomology of a symmetric product of an algebraic curve.
\newblock \textsl{Duke Math. J.} \textbf{108} (2001) 329--362

\bibitem[Br1]{Br} \textsc{S. Bradlow}: Vortices in holomorphic line
bundles over closed K\"ahler manifolds. \newblock \textsl{Commun. Math. Phys.}
\textbf{135} (1990), 1--17

\bibitem[Br2]{Br2} \textsc{S. Bradlow}: Special metrics and stability for
holomorphic bundles with global sections.
\newblock \textsl{J. Diff. Geom.}
\textbf{33} (1991), 169--213

\bibitem[BD1]{BD1} \textsc{S. Bradlow  {\upshape and} G. Daskalopoulos}: Moduli of stable pairs for holomorphic bundles over Riemann surfaces.
\newblock \textsl{Internat. J. Math.}
\textbf{2} (1991), 477--513

\bibitem[BD2]{BD2} \textsc{S. Bradlow  {\upshape and} G. Daskalopoulos}: Moduli of stable pairs for holomorphic bundles over Riemann surfaces II.
\newblock \textsl{Internat. J. Math.}
\textbf{4} (1993), 903--925

\bibitem[BDGW]{BDGW} \textsc{S. Bradlow, G.D. Daskalopoulos,  O. Garc\'{\i}a-Prada {\upshape and} R. Wentworth}: Stable 
augmented bundles over Riemann surfaces. In: N.J. Hitchin, P.E. Newstead, W.M. Oxbury (Eds.): {\it Vector Bundles in Algebraic Geometry}, LMS Lecture Notes
Series 208, Cambridge University Press, 1995; pp.\ 15--67

\bibitem[CCL]{CCL} \textsc{S.S. Chern, W.H. Chen {\upshape and} K.S. Lam}: 
\textit{Lectures on Differential Geometry}, World Scientific, 1999

\bibitem[D]{DonNS} \textsc{S.K. Donaldson}: A new proof of a theorem
of Narasimhan and Seshadri. \newblock \textsl{J. Diff. Geom.}
\textbf{18} (1983), 269--277

\bibitem[DK]{DonKro} \textsc{S.K. Donaldson {\upshape and} P. Kronheimer}: 
\textit{The Geometry of Four-Manifolds}, Clarendon Press, 1990

\bibitem[EH]{EisHar} \textsc{D. Eisenbud {\upshape and} J. Harris}: 
\textit{The Geometry of Schemes}, Springer, 2007

\bibitem[EINOS]{EINOS}
\textsc{M. Eto, Y. Isozumi, M. Nitta, K. Ohashi {\upshape and} N. Sakai}: 
Solitons in the Higgs phase: the moduli matrix approach.
\newblock \textsl{J. Phys. A: Math. Gen.} \textbf{39} (2006) R315--R392

\bibitem[F]{Ful} \textsc{W.\ Fulton}: 
\textit{Algebraic Topology: A First Course}, Springer, 1997

\bibitem[Ga]{GaDEP} \textsc{O. Garc\'{\i}a-Prada}: A direct existence proof for the
vortex equations over a compact Riemann surface.
\newblock \textsl{Bull. London Math. Soc.} \textbf{26} (1994), 88--96

\bibitem[GH]{GH} \textsc{P. Griffiths {\upshape and} J. Harris}: 
\textit{Principles of Algebraic Geometry},
John Wiley and Sons, 1978

\bibitem[Gr]{Gr} \textsc{A. Grothendieck}: Techniques de construction et th\'eor\`emes d'existence en g\'eom\'etrie alg\'ebrique. IV. Les sch\'emas de Hilbert,
S\'eminaire Bourbaki, Vol. 6, Exp. No. 221, 249--276, Soc. Math. France, Paris, 1995

\bibitem[HL]{HL} \textsc{D. Huybrechts {\upshape and} M. Lehn}: \textit{The Geometry
of Moduli Spaces of Sheaves}, Aspects of Mathematics, Vieweg \& Sohn, 1997

\bibitem[JT]{JafTau} \textsc{A. Jaffe {\upshape and} C. Taubes}: \textit{Vortices
and Monopoles}, Birkh\"auser, 1980

\bibitem[M]{McD}
\textsc{I.\ G.\ Macdonald}: Symmetric products of an algebraic curve.
\newblock \textsl{Topology} \textbf{1} (1962) 319--343

\bibitem[MN]{ManNasVVMS}
\textsc{N.S. Manton {\upshape and} S.M. Nasir}: Volume of vortex  
moduli spaces.
\newblock \textsl{Commun. Math. Phys.} \textbf{199} (1999) 591--604

\bibitem[MR]{ManRom}
\textsc{N.S. Manton {\upshape and} N.M. Rom\~ao}: Vortices and Jacobian varieties.
\newblock {\tt arXiv:1010.0644}

\bibitem[OT]{OkoTel}
\textsc{Ch. Okonek {\upshape and} A. Teleman}: Gauge theoretical equivariant
Gromov--Witten invariants  and the full Seiberg--Witten invariants of ruled surfaces.
\newblock \textsl{Commun. Math. Phys.} \textbf{227} (2002) 551--585

\bibitem[Th]{Tha}
\textsc{M. Thaddeus}: Stable pairs, linear systems, and the Verlinde formula.
\newblock \textsl{Invent. Math.} \textbf{117} (1994) 317--353

\bibitem[To]{TonQVS}
\textsc{D. Tong}: Quantum vortex strings: a review.
\newblock \textsl{Ann. Phys.} \textbf{324} (2009) 30--52

\bibitem[UY]{UhlYau}
\textsc{K. Uhlenbeck {\upshape and} S.T. Yau}: On the existence of Hermitian-Yang-Mills
connections in stable vector bundles.
\newblock \textsl{Comm. Pure Appl. Math.} \textbf{39} (1986) S257--S293

\end{small}

\end{thebibliography}

\end{document}